\documentclass[11pt]{article} 

\usepackage[margin=2.8cm]{geometry}
\usepackage{setspace}

\usepackage{color}
 \usepackage{amsfonts,amssymb,latexsym,bm,amsmath,amsfonts}
 \usepackage[latin1]{inputenc}  

\newenvironment{thmbis}[1]
  {\renewcommand{\thetheorem}{\ref{#1}${}^*$}%
   \addtocounter{theorem}{-1}%
   \begin{theorem}}
  {\end{theorem}}

\newtheorem{theorem}{Theorem}[section]

\newtheorem{lemma}{Lemma}[section]

\def\beq{\begin{equation}}  \def\eeq{\end{equation}}
\def\bb{\begin{eqnarray*}}  \def\ee{\end{eqnarray*}}
\def\b{\begin{eqnarray}}    \def\e{\end{eqnarray}}
\def\1{\hbox{\rm\setbox1=\hbox{1}\copy1\kern-.5\wd1 I}}
\def\D{{\hbox{\rm\setbox1=\hbox{I}\copy1\kern-.45\wd1 D}}}

\def\R{\mathbb{R}}
\def\N{\hbox{\rm\setbox1=\hbox{I}\copy1\kern-.55\wd1 N}}
  
\def\L{{\cal L}} 

\def\p{\hbox{\rm\setbox1=\hbox{I}\copy1\kern-.45\wd1 P}}
\def\<{{^{_<}}}   \def\>{{^{_>}}}   
     
   \def\bib{\vspace{-2mm}\bibitem}

\def\l{\lambda}

\def\({\left(}  
\def\){\right)}

\newcommand{\eec}{{\mathrm e}} 	\newcommand{\RR}{\mathbb{R}}
\newcommand{\Exponent}[1]{\exp\Biggl\{#1\Biggr\}} 

\newcommand{\dtv}{d_{_{TV}}}
\newcommand{\dloc}{d_{_{0}}}    					

\newcommand{\dw}{d_{_{W}}}                           

\newcommand{\norm}[1]{\|#1\|}               		

\newcommand{\tvnorm}[1]{\norm{#1}_{_{TV}}}

\newcommand{\lnorm}[1]{\norm{#1}_{_\mathrm{0}}}  
\newcommand{\dlrnorm}[1]{\norm{#1}_r}

\newcommand{\wnorm}[1]{\ab{#1}_{_{W}}} 	
\newcommand{\Lnorm}[1]{\ab{#1}_r}

\newcommand{\ab}[1]{\vert#1\vert}           
\newcommand{\ii}{{\mathrm i}}               
\newcommand{\dd}{{\mathrm{d}}}              
		
\newcommand{\dirac}{I}                      
\newcommand{\F}{{\mathcal{F}}} 							
\newcommand{\ZZ}{\mathbb{Z}}	
	\newcommand{\w}{\widehat}

\newcommand{\eit}{\eec^{\ii t}}							

\newcommand{\M}{\cal M}

\newcommand{\A}{\alpha}
\newcommand{\B}{\beta}

\begin{document}

\title{Skellam compound Poisson approximation to the sums of symmetric Markov dependent random variables}

\author{ V. \v Cekanavi\v cius\footnote{Corresponding author.}  \, and  G. Liaudanskait\.e \\
 {\small
Institute of Applied Mathematics, Faculty  of Mathematics and Informatics, Vilnius University,}\\
{\small Naugarduko 24, Vilnius 03225, Lithuania.}\\{\small E-mail:
vydas.cekanavicius@mif.vu.lt } {\small and  gabija.liaudanskaite@mif.vu.lt}
 }
\date{}

\maketitle

\begin{abstract}

The sum of symmetric Markov dependent three-point random variables is approximated by the difference of two independent Poisson random variables (Skellam random variable).  The accuracy is estimated in local, total variation and Wasserstein metrics. Properties of convolutions of measures is used for the proof.

\vspace*{.5cm} \noindent {\emph{Key Words:} \small  compound Poisson approximation; concentration function; local metric; Markov chain, Skellam distribution; total variation metric; Wasserstein metric}

\vspace*{.5cm} \noindent {\small {\bf MSC 2000 Subject
Classification}:
Primary 60F05.   
Secondary  60J10; 
}
\end{abstract}

\newpage

It is known that symmetry of random variables (r.v.s) can radically improve compound Poisson approximation for sums of  independent r.v., see, for example, \cite{Pre85}, Arak's inequality  \cite{A80} or its generalization to higher dimensions  \cite{GZ21,Z03}. The case of dependent r.v.s is less explored. In this paper, we use the difference of two independent Poisson variables (Skellam r.v.) to approximate  symmetric three state Markov chain.  Under quite general assumptions (for example, for the classical setting of sequences of r.v.s) the order of approximation is $O(n^{-1})$. The accuracy is improved to $O(n^{-2})$ by the first order asymptotic expansion.

Let	 $\R$ denote the set of all real numbers, $\,\ZZ\,$ denote the set of all integers, 
 $\N$ denote the set of natural numbers. 
 The set of all one-dimensional distributions  is denoted by $\F$. Notation $\F_\ZZ$ is used for distributions 
 concentrated on $\ZZ$.          , 
 Above notation can be extended in a natural way to $\M$  
  for bounded signed measures, i.e. measures  such that
$ -\infty<\inf_{A}M\{A\}\leq\sup_{A}M\{A\}<\infty$,
where infimum and supremum are taken over all Borel sets, and $M\{A\}<0$ is allowed.   Measures concentrated on $\ZZ$ are denoted by $\M_{\ZZ}\subset\M$.

Let $X$ be random variable (r.v.). We denote  distribution of $X$ by $\L(X)$.
If $F=\L(X)$, then for any Borel set $A$, $F\{A\}=\p(X\in A)$. To make notation shorter we write $F\{k\}:=F\{\{k\}\}$. We denote by $\pi_\l\sim{\cal P}(\l)$  Poisson r.v. with parameter $\l>0$. The difference $\pi_{\l_1}-\tilde\pi_{\l_2}$ of two independent Poisson r.v.s. $\pi_{\l_1}$ and $\tilde\pi_{\l_2}$ is called Skellam r.v., see \cite{GK18,Sk46}. Note that probabilities of Skellam r.v. can be expressed through modified Bessel function of the first  kind:
\begin{equation}\label{SkP}
\p(\pi_{\l_1}\!-\!\tilde\pi_{\l_2}=k)=\eec^{-\l_1-\l_2}(\l_1/\l_2)^{k/2}{\rm I}_{\ab{k}}(2\sqrt{\l_1\L_2}), \quad k\in\ZZ.\end{equation} 
For other expressions of Skellam distribution and some properties of Bessel functions, see \cite{AO10}.
Due to \cite{I37} Skellam distribution with $\l=\tilde\l$ is sometimes called Irwin distribution.
 
Skellam distribution is one of the simplest cases of compound Poisson  distributions. In general compound Poisson (CP) distribution 
is distribution of a random sum of independent identically distributed (iid)  r.v.s, when number of summands is Poisson r.v. More precisely, r.v.
 $\,Y\,$ has a CP distribution
if
 \beq                                         \label{CP} 
 Y \stackrel{d}{=} \sum_{k=0}^{\pi_\l}X_i\,,
 \eeq
where Poisson  random variable $\pi_\l$ is independent of
$\{X_i\}_{i\ge1}$, $X_0\!=\!0,$ and random variables $X_1,X_2,...$ are iid r.v.s., i.e. $X_i\stackrel{d}{=}X\,$ $\,(i\!\ge\!1)$.
Distribution $\L(X)$ is called \emph{compounding distribution.} 

Arguably the most popular metric used to estimate closeness of lattice distributions is the total variation metric 
$\dtv(X;Y)=\sup_A\ab{\p(X\!\in\!A)-\p(Y\!\in\!A)}$, where supremum is taken over all Borel sets.
Other popular metrics include 
point (local) metric
$\dloc(X;Y)= \sup_{k\in\ZZ} \ab{\p(X\!=\!k)-\p(Y\!=\!k)}$ and Wasserstein (a.k.a. Gini, a.k.a Kantorowich, a.k.a. Fortet-Mourier, a.k.a. $L_1$) metric
$\dw(X;Y)=\sum_{k=-\infty}^\infty\ab{\p(X\!\leq\!x)-\p(Y\!\leq\!x)}$. For equivalent definitions involving supremum over special class of functions, see \cite{BHJ}, App. A.1 or \cite{CN22}, Sec. 1.2. 

 In this paper we use  properties of measures and for us is more convenient to use norm notation.   For $M\in\M_\ZZ$  we define  local, total variation, and Wasserstein norms respectively by:
  \[
 \lnorm{M}:=\sup_{k\in\ZZ}\ab{M\{k\}},\quad \tvnorm{M}:=\sum_{k=-\infty}^\infty\ab{M\{k\}},\quad 
\wnorm{M}:=\sum_{k=-\infty}^\infty\ab{M\{(-\infty,k]\}}.
  \]

\noindent Relations between norms and corresponding metrics are quite direct. Let $\L(X),\L(Y)\in\F_\ZZ$, then
 \[
 \dtv(X;Y)=\frac{1}{2}\tvnorm{\L(X)-\L(Y)},\quad \dloc(X;Y)=\lnorm{\L(X)-\L(Y)},\quad
  \dw(X;Y)=\wnorm{\L(X)-\L(Y)}.
 \]
 Note that total variation metric is {half} of the total variation norm. 
Total variation norm and Wasserstein norms can be extended  respectively  to
$\ell_r$ and $L_r$, $(r\geq 1)$ norms:
\[\dlrnorm{M}:=\left(\sum_{k=-\infty}^\infty\ab{M\{k\}}^r\right)^{1/r},
 \quad
 \Lnorm{M}:=\left(\sum_{k=-\infty}^\infty \ab{M\{(-\infty,k]\}}^r\right)^{1/r}.\]

If $M,V\in \M_\ZZ$, then  $M*V\in\M_\ZZ$ denotes their convolution: for any $m\in\ZZ$
\[M\!*\!V\{m\}=\sum_{k=-\infty}^\infty M\{m\!-\!k\}V\{k\}.\]
Powers of measures are understood in convolution sense $M^{*k}=M\!*\!M\!*\!\cdots\!*\!M$.
For any two independent r.v.s $X$, $Y$ the distribution of their sum equals convolution of corresponding distributions: $\L(X\!+\!Y)=\L(X)\!*\!\L(Y)$.

Let $I_a$ denote the distribution  of degenerate r.v. $X$ which is concentrated at $a\in\R$, i.e. $\p(X\!=\!a)=1$, or in terms of measures,
	$I_a\{a\}=1$, $I\{\R\setminus\{a\}\}=0$.
Further on $I\equiv I_0$. Observe that $(I_1)^{*k}=I_k$.  Assuming  $M^{*0}\equiv I$ we can define exponential measure 
	\[\exp\{M\}=\sum_{k=0}^\infty \frac{M^{*k}}{k!}.\]
Note that any CP distribution can be written as exponential measure. For example, distribution of r.v. $Y$ defined in (\ref{CP}) can be written as
$\L(Y)=\exp\{\l(\L(X)-I)\}$. Poisson distribution and Skellam distribution can resp. be written as
\[\L(\pi_\l)=\exp\{\l(I_1\!-\!I)\},\quad \L(\pi_{\l_1}-\tilde\pi_{\l_{2}})=\exp\{\l_1(I_1\!-\!I)+\l_{2}(I_{-1}\!-\!I)\}.\]

Fourier transform of $M\in\M_\ZZ$ is defined as
\[\w M(t)=\sum_{k=-\infty}^\infty \eec^{\ii t k}M\{k\}, \qquad t\in \R.\]
Here $\ii$ denotes imaginary unit $\ii^2=-1$. If $F\in\F$, then $\w F(t)$ is its characteristic function. Observe that characteristic function of Poisson r.v. $\pi_\l\sim{\cal P}(\l)$ is equal to
$\exp\{\l(\eit\!-\!1)\}$ and characteristic function of CP r.v. defined in (\ref{CP}) is equal to $\exp\{\lambda(\w F(t)-1)\}$.

We use the same symbol $C$ to denote different absolute constants.  
 Symbol $\Theta$  denotes any measure with total variation $\tvnorm{\Theta}\!\leq\!1$. 


\section{Known results}

Presman \cite{Pre85} considered Skellam approximation to the sum of independent symmetrized Bernoulli variables.
For the sum of independent three point r.v.s Presman's result can be reformulated
   in the following way:
Let $X,X_1,X_2,\dots,X_n$ be iid r.v.s, $\p(X=-1)=\p(X_i=1)=p\leq 1/4$, $\p(X=0)=1\!-\!2p$, $\tilde S_n=X_1+\cdots+X_n$ and let $\pi_{np}$,$ \tilde \pi_{np}$ be two independent Poisson r.v.s; independent of $X_i$. Then
\begin{equation}\label{sympre83}
\tvnorm{\L(\tilde S_n)-\L(\pi_{np}\!-\!\tilde\pi_{np})}\leq C\min(np^2; n^{-1}). 
\end{equation}
 The accuracy of approximation in (\ref{sympre83}) is of correct order.  Estimate (\ref{sympre83}) is closely related to the so-called first uniform Kolmogorov theorem on the best possible infinitely divisible approximation to the sum of independent r.v.s , see \cite{AZ} Introduction and Chapter IV. 

Natural generalization of Binomial distribution is the sum of Markov dependent Bernoulli variables (Markov Binomial r.v.). Poisson and CP approximations to Markov Binomial distribution are well investigated  \cite{CeVe10,OSV08,XZ09}. On the other hand, there are just a few results for CP approximation of symmetric  Markov dependent r.v.s. In \cite{SC16} Markov dependent analogue  of $\tilde S_n$ was considered.
Let $\xi_0, \xi_1, \dots \xi_n, \dots$ be a non-stationary three state $\{ a_1, a_2, a_3 \}$ Markov chain. We denote the distribution of $S_n=f(\xi_1)+\dots+f(\xi_n)$ $(n\in\N)$,  where $f(a_1)=-1$, $f(a_2)=0$, $f(a_3)=1$, by $F_n=\L(S_n)$. The initial distribution is $\p(\xi_0=a_1)=p_1$, $\p(\xi_0=a_2)=p_2$ and $\p(\xi_0=a_3)=p_3$. For $\A,\B\in (0;0.5)$ we define the following transition  probabilities:
\begin{eqnarray*}
&\displaystyle \p(\xi_i=a_1\,|\,\xi_{i-1}=a_1)=\A,\quad
\p(\xi_i=a_2\,|\,\xi_{i-1}=a_1)=1-2\A,\quad
\p(\xi_i=a_{3}\,|\,\xi_{i-1}=a_1)=\A,&\\
&\displaystyle \p(\xi_i=a_1\,|\,\xi_{i-1}=a_2)=\B,\quad
\p(\xi_i=a_2\,|\,\xi_{i-1}=a_2)=1-2\B,\quad
\p(\xi_i=a_{3}\,|\,\xi_{i-1}=a_2)=\B,&\\
&\displaystyle \p(\xi_i=a_1\,|\,\xi_{i-1}=a_3)=\A,\quad
\p(\xi_i=a_2\,|\,\xi_{i-1}=a_3)=1-2\A,\quad
\p(\xi_i=a_{3}\,|\,\xi_{i-1}=a_3)=\A.&
\end{eqnarray*}
Note that if $\A=\B=p$, then $S_n$ becomes the sum of independent three point r.v.s  from (\ref{sympre83}), $S_n\stackrel{d}{=}\tilde S_n$. 

Next we formulate the main result from \cite{SC16}, where it was assumed that
\begin{equation}
0\leq \A \leq \frac{1}{30},\quad 0\leq \B \leq \frac{1}{30}.\label{ab}
\end{equation}   
and the following compound distributions introduced:
\begin{eqnarray*}
L&=&\frac{1}{2}(\dirac_{-1}+\dirac_{1}),\quad H=(1-2\A)L*\displaystyle\sum_{j=0}^{\infty}(2\A L)^{*j},
\quad
E=\bigg( 1-\frac{2\A p_2}{1-2\A}\bigg)\dirac+\frac{2\A p_2}{1-2\A}L,\\
G&=&\Exponent{\frac{2\B(1-2\A)}{1-2\A+2\B}(H-\dirac)},  
\quad
K=\frac{1-2(\A-\B)}{1+2\B}\displaystyle\sum_{j=0}^{\infty}\Big( \frac{2\A}{1+2\B}L \Big)^{*j}.  
\end{eqnarray*}
%
%

\begin{theorem} \label{T1} (\cite{SC16}) Let condition (\ref{ab}) hold. Then, for all $n=1,2,\dots,$
\begin{eqnarray*}
\tvnorm{F_n-E\!*\!K\!*\!G^{*n}}&\leq& C\big(\min(n^{-1};\B)+0.2^n\ab{\A-\B}\big),\label{1.1}
\\
\lnorm{F_n-E\!*\!K\!*\!G^{*n}}&\leq& C\big( \min(n^{-3/2}\B^{-1/2};\B  )+0.2^n\ab{\A-\B}\big),\label{1.12}\nonumber
\\
\wnorm{F_n-E\!*\!K\!*\!G^{*n}}&\leq& C\big( \min(n^{-1/2}\B^{1/2};\B)+0.2^n\ab{\A-\B}\big).\nonumber
\end{eqnarray*}
\end{theorem}
Though the accuracy of estimates in Theorem \ref{T1} is comparable to (\ref{sympre83}), approximation $E*M*G^{*n}$ is structurally very complicated, since $G$ is CP distribution with compounding compound geometric distribution. More precisely $E*M*G^{*n}=\L(Y_1\!+\!Y_2\!+\!Y_3)$, where $Y_1,Y_2,Y_3$ are independent r.v.s and
$\p(Y_1\!=\!1)=\p(Y_1\!=\!-1)=\A p_2/(1\!-\!2\A)$, $\p(Y_1\!=\!0)=1-2\A p_2/(1\!-\!2\A)$; $Y_2=\sum_{k=0}^N\eta_k$, where $\eta;\eta_1,\eta_2,\dots$ are iid r.v.s, $\p(\eta=-1)=\p(\eta=1)=1/2$ and $N$ is geometric r.v. independent of $\eta_i$:
$\p(N\!=\!k)=\big({2\A}/(1\!+\!2\B)\big)^k\big(1-2\A/(1\!+\!2\B)\big)$, $k=0,1,2\dots$
 We have $\L(Y_3)=G^{*n}$, therefore 
 \[Y_3=\sum_{j=0}^{\pi_{\tilde\l}}Y_{3j},\quad\tilde\l=2n\B(1\!-\!2\A)/(1\!-\!2\A\!+\!2\B)\]
  and $Y_{30},Y_{31},Y_{32},\dots$ are iid r.v.s independent of $\pi_{\tilde\l}$. Moreover, $Y_{30}$ is also a random sum with geometric number of summands:
  \[Y_{30}=\sum_{m=1}^{\tilde N}\eta_m,\quad \p(\tilde N=k)=(1\!-\!2\A)(2\A)^{k-1},\quad k=1,2,\dots\]
Here all r.v.s are independent.
Such complicated structure makes approximation $E*M*G^{*n}$ of limited practical use. In \cite{SC16}, first order approximation with even more complicated  structure was also used. 

Theorem \ref{T1} is natural for small values of  $\B$ such as $\B=O(1/n)$ and $\A=O(1)$. On the other hand, $E*M*G^{*n}$ is structurally overcomplicated if $\B$ is of similar magnitude to $\A$. For example, if $\A=\B$, then $F_n=\L(\tilde S_n)$ from  (\ref{sympre83}), but  $G^{*n}$ still remains CP with compounding compound geometric distribution. Therefore  Theorem \ref{T1} is not a direct extension of 
estimate (\ref{sympre83}) to Markov chain.

\section{Results}

Our goal is to apply Skellam approximation to the sum of symmetric Markov dependent  r.v.s  defined in previous Section. Let
$\l=\B/(1-2\A+2\B)$ and denote symmetric Skellam distribution by 
\[D=\L(\pi_\l\!-\!\tilde\pi_\l)=\exp\bigg\{\frac{\B}{1-2\A+2\B}(I_1\!-\!I+I_{-1}\!-\!I)\bigg\}.\]
Here $\pi_\l,\tilde\pi_\l$ are two independent Poisson r.v.s with parameter $\l$.

\begin{theorem}\label{T1GC}  Let condition (\ref{ab}) hold 
Then
\begin{eqnarray*}
\tvnorm{ F_n- D^{*n}}&\leq& \frac{C}{n}\bigg(1+\frac{\ab{\A-\B}}{\B}\bigg),
\\
\lnorm{ F_n- D^{*n}}&\leq& \frac{C}{n\sqrt{n\B}}\bigg(1+\frac{\ab{\A-\B}}{\B}\bigg),
\\
\wnorm{ F_n- D^{*n}}&\leq& C\sqrt{\frac{\B}{n}}\bigg(1+\frac{\ab{\A-\B}}{\B}\bigg).
\end{eqnarray*}
\end{theorem}

If $\A=\B$ then $D$ coincides with Skellam approximation used in (\ref{sympre83}). Thus Theorem \ref{T1GC} is direct extension of (\ref{sympre83}) to Markov dependent r.v.s for $\B\geqslant 1/n$. In many cases (such as $\A/\B\leq C$) the estimates in Theorem \ref{T1GC} are of the same order as estimates of Theorem \ref{T1}.

The accuracy of approximation can be improved by first order asymptotic expansion. Let

\begin{eqnarray*}
A_1&:=&\frac{2(\A\!-\!\B)}{1\!-\!2\A\!+\!2\B}\bigg(\frac{1\!-\!2\A}{1\!-\!2\A\!+\!2\B}-p_2\bigg)(L\!-\!I),\\
A_2&:=&\frac{2\B}{(1\!-\!2\A\!+\!2\B)^2}\bigg(
\frac{2(\A\!-\!\B)(1\!-\!2\A)}{1\!-\!2\A\!+\!2\B}-\B\bigg)(L\!-\!I)^{*2}.
\end{eqnarray*}
We recall that $2(L\!-\!I)=I_1\!-\!I+I_{-1}\!-\!I$.

\begin{theorem}\label{T2} Let condition (\ref{ab}) hold. 
Then
\begin{eqnarray*}
\tvnorm{ F_n- D^{*n}*(I+A_1+nA_2)}&\leq& \frac{C}{n^2}\bigg(1+\frac{(\A-\B)^2}{\B^2}\bigg),
\\
\lnorm{ F_n- D^{*n}*(I+A_1+nA_2)}&\leq& \frac{C}{n^2\sqrt{n\B}}\bigg(1+\frac{(\A-\B)^2}{\B^2}\bigg),
\\
\wnorm{ F_n- D^{*n}*(I+A_1+nA_2)}&\leq& \frac{C}{n}\sqrt{\frac{\B}{n}}\bigg(1+\frac{(\A-\B)^2}{\B^2}\bigg).
\end{eqnarray*}
\end{theorem}
Easily verifiable identities for any $M\in\M_\ZZ$, $k\in\ZZ$:
\begin{equation}\label{dif}
(I_1-I)*M\{k\}=M\{k\!-\!1\}-M\{k\}, \quad (I_{-1}-I)*M\{k\}=M\{k\!+\!1\}-M\{k\}\end{equation}
allow to express 'probabilities' of asymptotic expansion from Theorem \ref{T2} through forward and backward differences of probabilities of Skellam distribution $D$, for which relation (\ref{SkP}) can be used.

The method of proof does not allow to get small constants. On the other hand, the idea of their magnitude can be understood if we combine Theorem \ref{T2} with Theorem \ref{T1GC}.

\begin{thmbis}{T1GC}\label{T1star} Let condition (\ref{ab}) hold. 
Then
\begin{eqnarray*}
\tvnorm{ F_n- D^{*n}}&\leq&   \frac{0.61}{n}\bigg(1+\frac{3.21\ab{\A-\B}}{\B}\bigg)+\frac{C}{n^2}\bigg(1+\frac{(\A-\B)^2}{\B^2}\bigg),
\\
\lnorm{ F_n- D^{*n}}&\leq& \frac{0.6}{n\sqrt{n\B}}\bigg(1+\frac{3\ab{\A-\B}}{\B}\bigg)+\frac{C}{n^2\sqrt{n\B}}\bigg(1+\frac{(\A-\B)^2}{\B^2}\bigg),
\\
\wnorm{ F_n- D^{*n}}&\leq& 0.5\sqrt{\frac{\B}{n}}\bigg(1+\frac{3.9\ab{\A-\B}}{\B}\bigg)+\frac{C}{n}\sqrt{\frac{\B}{n}}\bigg(1+\frac{(\A-\B)^2}{\B^2}\bigg).
\end{eqnarray*}
\end{thmbis}

Estimates for total variation and Wasserstein norms can be extended to estimates for $\ell_r$ and $L_r$ norms. 

\begin{theorem}\label{lrLr} Let condition (\ref{ab}) hold and let $r\geq 1$. 
Then
\begin{eqnarray*}
\dlrnorm{F_n- D^{*n}}&\leq& Cn^{-(3r-1)/2r}\B^{-(r-1)/2r}(1+{\ab{\A-\B}}/{\B}),\\
  \Lnorm{F_n-D^{*n}}&\leq& Cn^{(2r-1)/2r}\B^{1/2r}(1+\ab{\A-\B}/\B).
\end{eqnarray*}
\end{theorem}
The remaining part of the paper is devoted to the proofs of theorems.     . 

\section{Auxiliary results}

  For $M,V\in\M_\ZZ$, the following relations hold
\begin{eqnarray}
&&\lnorm{M}\leq \tvnorm{M},\quad \tvnorm{M*V}\leq\tvnorm{M}\tvnorm{V}, \quad  \tvnorm{M}\leq\exp\{\tvnorm{M}\}, \label{baznrm}\\
&&\wnorm{M*V}\leq \wnorm{M}\tvnorm{V},\label{nrmrel}\\
&&  \lnorm{M}\leq\tvnorm{(I_1\!-\!I)\!*\!M},\quad \tvnorm{M}=\wnorm{(I_1\!-\!I)*M}. 
\label{locWToTV}
 \end{eqnarray}
Estimates (\ref{baznrm}) are well-known an hold  for even more general $M,V\in\M$. For the proof of   (\ref{nrmrel}) observe that
\begin{eqnarray*}
\wnorm{M*V}&:=&\sum_{k\in\ZZ}\Big\vert\sum_{j\in\ZZ}M\{(-\infty,k\!-\!j]\}V\{j\}\Big\vert\leq
\sum_{k\in\ZZ}\sum_{j\in\ZZ}\ab{M\{(-\infty,k\!-\!j]\}}\ab{V\{j\}}\nonumber\\
&=&\sum_{j\in\ZZ}\ab{V\{j\}}\sum_{k\in\ZZ}\ab{M\{(-\infty,k\!-\!j]\}}=\wnorm{M}\tvnorm{V}. \label{WTV}
\end{eqnarray*}
First estimate in (\ref{locWToTV}) is Eq. (3.23) in \cite{C16}. Second estimate follows directly from definition of the  Wasserstein norm and (\ref{dif}).

Next we present some facts about exponential measures and CP distributions. For any numbers $u_1,u_2\in\RR$, and $F\in\F$ exponential measures satisfy simple relation
\[\exp\{(u_1+u_2)(F-I)\}=\exp\{u_1(F-I)\}*\exp\{u_2(F-I)\}.\]
We recall that to symbol $\Theta$ is used for any measure satisfying $\tvnorm{\Theta}\leq 1$ . Thus,  if $F\in\F$, $u_1,u_2>0$, then we can write, for example, the following equality
\[\exp\{(u_1+u_2)(F\!-\!I)\}=\exp\{u_1(F\!-\!I)\}*\Theta,\] 
since 
\[ \tvnorm{(I_1-I)*\exp\{u_2(I_1-I)\}}\leq 2\tvnorm{\exp\{u_2(I_1-I)\}}= 2.\]
Here we  used the fact that $\exp\{u_2(I_1-I)\}\in\F$  and, therefore, its total variation is equal to 1.

 The next lemma gives some properties of exponential measures that can be easily verified by integration.


\begin{lemma}\label{bazine} Let $V,M\in\M$, , $k\in\N$. Then
\[
\eec^{M}-\eec^{V}=(M\!-\!V)*\int_0^1\eec^{\tau M+(1-\tau)V}\dd\tau,\quad\eec^V=I+\sum_{j=1}^k\frac{V^{*j}}{j!}+\frac{V^{*(k+1)}}{k!}*\int_0^1\eec^{\tau V}(1-\tau)^k\dd\tau.
\]
\end{lemma}
 
 The following lemma was proved in \cite{C16}, p. 31.
 

 \begin{lemma} (\cite{C16}) \label{L2.5} Let $\lambda>0$, $F\in\F$. Then
 \[\bigg\|\exp\bigg\{\lambda (F\!-\!I)+\frac{2\lambda}{7}(F\!-\!I)^{*2}*\Theta\bigg\}\bigg\|\leq C.\]
  \end{lemma}
 

 \begin{lemma} \label{Roslema} Let $t>0$, $F\in\F$, $k\in N$. Then
 \begin{eqnarray}
&&\tvnorm{(F\!-\!I)*\exp\{t(F-I)\}}\leq\sqrt{2/\eec t},\quad 
\tvnorm{(F\!-\!I)^{*2}*\exp\{t(F-I)\}}\leq 3/\eec t,\label{C.2}\\  
&&\tvnorm{(F\!-\!I)^{*k}*\exp\{t(F-I)\}}\leq \sqrt{k!}\,t^{-k/2}.\label{C.3}
 \end{eqnarray} 
 \end{lemma}
 Inequalities (\ref{C.2})  are proved  in \cite{Roos01}, Lemma 3. Estimate (\ref{C.3}) is from Lemma 4, \cite{R2003}.

For estimation of local constants we will use Lemma 4.6 from \cite{CeRoAISM06}.
\begin{lemma} (\cite{CeRoAISM06}) \label{zloc} Let $j=1,2,\dots$ and $t\in(0,\infty)$. If $F$ is symmetric distribution concentrated on the set $\ZZ\setminus\{0\}$, then
\[
\lnorm{(F-I)*\exp\{t(F-I)\}}\leq 2\bigg(\frac{j+1/2}{t\eec}\bigg)^{j+1/2}. \]
\end{lemma}


 In \cite{BRG51} the so-called Bergstr\"om's identity was proved.
 
 \begin{lemma} \label{bergidentL} (\cite{BRG51}) Let $V,M\in\M$, $k,n\in\N$. Then
 \begin{eqnarray*}
 V^{*n}-M^{*n}\!&=&\!\sum_{m=1}^k\binom{n}{m}(V-M)^{*m}*M^{*(n-m)}\\
 \!&+&\!(V-M)^{*(k+1)}*\sum_{m=k}^{n-1}\binom{m}{k}V^{*(n-1-m)}*M^{*(m-k)}.
 \end{eqnarray*}
 \end{lemma}
 Bergstr\"om's identity is usually combined with combinatorial identity
 \begin{equation*}\label{bin}
 \sum_{m=k}^{n-1}\binom{m}{k}=\binom{n}{k+1}.
 \end{equation*}


Next we present some results from \cite{SC16}.  We recall that $L=\frac{1}{2}(\dirac_{-1}+\dirac_{1})$ and set
$U=L-I$.

\begin{lemma}\label{3Snotation} (\cite{SC16}) Let (\ref{ab}) be satisfied. Then
\[{F}_n=P_1*\Lambda_1^{*n}* W_1+ P_2*\Lambda_2^{*n}*W_2\label{F-n1},\]
 where
\begin{eqnarray*}
P_1\!\!&=&\!\!\frac{\pi_2}{1\!-\!2\A}(\Lambda_1\!-\!I-2\A U),\qquad
P_2=\frac{\pi_2}{1\!-\!2\A}(\Lambda_2\!-\!I-2\A U),
\\
\Lambda_{1,2}\!\!&=&\!\!\frac{1}{2}\bigg( (1\!+\!2\A\!-\!2\B)I+2\A U \pm((1\!-\!2\A\!+\!2\B)I -2\A U)*
\displaystyle\sum_{j=0}^{\infty}\binom{1/2}{j} \Delta^{*j}  \bigg),
\\
W_{1,2}\!\!&=&\!\!\frac{1}{2}\bigg(\dirac\pm\bigg(I+\frac{2\A U}{1\!-\!2\A\!+\!2\B}\bigg)*
 \displaystyle\sum_{j=0}^{\infty}\Big( \frac{2\A}{1\!-\!2\A\!+\!2\B}\Big)^jU^{*j}
\displaystyle*\sum_{j=0}^{\infty}\binom{-1/2}{j} \Delta^{*j} \bigg),\\
\Delta\!\!&=&\!\!\frac{8\B U}{(1\!-\!2\A\!+\!\B)^2}*\bigg(\displaystyle\sum_{j=0}^{\infty}\bigg( \frac{2\A}{1\!-\!2\A\!+\!2\B}\bigg)^{*j}U^{*j}\bigg)^2,\quad \tvnorm{\Delta}\leq 0.62,\\ 
\tvnorm{\Lambda_2}^n\!\!&\leq&\!\! 15.5\ab{\A\!-\!B}0.2^n,\quad\tvnorm{W_i}\leq C,\quad\wnorm{W_i}\leq C, \quad\tvnorm{P_i}\leq C,\ i=1,2. 
\end{eqnarray*}
\end{lemma}
Note that  expressions of $W_i$ and $\Delta$ in above are based on their Fourier transforms from \cite{SC16}, p. 421.


Lemma \ref{3Snotation} allows to write expansions in convolution powers of $U$.

\begin{lemma}\label{3Smeasures} Let (\ref{ab}) be satisfied. Then
\begin{eqnarray}
U&=&\frac{1}{2}(I_1\!-\!I)^{*2}*I_{-1},\label{X2}\\
\Delta&=&{\displaystyle \frac{8\B U}{(1\!-\!2\A\!+\!2\B)^2}+\frac{32\A\B U^{*2}}{(1\!-\!2\A\!+\!2\B)^3}+C\A^2\B U^{*3}*\Theta,   }\label{Dlta}\\
\Lambda_1&=&{ I+\frac{2\B U}{1\!-\!2\A\!+\!2\B}+\frac{4\B(\A\!-\!\B)(1\!-\!2\A)U^{*2}}{(1\!-\!2\A\!+\!2\B)^3}+C(\B^3\!+\!\A^2\B\!+\!\A\B^2)U^{*3}*\Theta,}\label{Lmbda12}\\
P_1&=&I+\frac{2\pi_2(\B\!-\!\A)U}{1\!-\!2\A\!+\!2\B}+C(\B^2\!+\!\A\B)\pi_2U^{*2}*\Theta,\label{PP1}\\
W_1&=&I+\frac{2(\A\!-\!\B)(1\!-\!2\A)U}{(1\!-\!2\A\!+\!2\B)^2}+C(\A^2\!+\!\B^2)U^{*2}*\Theta,\label{WW1}\\
W_2&=&C(\A\!+\!\B)U*\Theta.\label{WW2}
\end{eqnarray}
\end{lemma}

Observe that form Lemma \ref{3Smeasures} shorter expansions immediately follow. For example,
\begin{equation}\label{Lambda1short}
\Lambda_1=I+\frac{2\B U}{1\!-\!2\A\!+\!2\B}+C(\B^2\!+\!\A\B)U^{*2}*\Theta.
\end{equation}

\textbf{Proof of Lemma \ref{3Smeasures}}. Expression (\ref{X2}) is simple identity. Total variation of any distribution is equal to 1. Therefore, $\tvnorm{U}\leq \tvnorm{L}+\tvnorm{I}=2$.  
From (\ref{ab}) it follows that
\[\frac{2\A}{1\!-\!2\A\!+\!2\B}\leq \frac{2\A}{1\!-\!2\A}\leq \frac{1}{7}.\]
Therefore
\[
\sum_{j=0}^{\infty}\bigg( \frac{2\A}{1\!-\!2\A\!+\!2\B}\bigg)^{*j}U^{*j}=I+\frac{2\A U}{1\!-\!2\A\!+\!2\B}+C\A^2U^{*2}*\Theta
\]
and
\[\bigg(I+\frac{2\A U}{1\!-\!2\A\!+\!2\B}+C\A^2U^{*2}*\Theta\bigg)^{*2}=I+\frac{4\A U}{1\!-\!2\A\!+\!2\B}+C\A^2U^{*2}*\Theta.
\]
Combining the last expression with Lemma \ref{3Snotation} we get (\ref{Dlta}). Observe that
\[
\sum_{j=3}^\infty\bigg\vert\binom{1/2}{j}\bigg\vert \tvnorm{\Delta}^{j-3}\leq\sum_{j=3}^\infty 0.62^{j-3}\leq 2.64;\quad \sum_{j=2}^\infty\bigg\vert\binom{-1/2}{j}\bigg\vert \tvnorm{\Delta}^{j-2}\leq 2.64.\] 
Therefore,
\[
\sum_{j=0}^{\infty}\binom{1/2}{j} \Delta^{*j}= I+\frac{1}{2}\Delta-\frac{1}{8}\Delta^{*2}+C\Delta^{*3}*\Theta,\quad
\sum_{j=0}^{\infty}\binom{-1/2}{j} \Delta^{*j}= I-\frac{1}{2}\Delta+C\Delta^2*\Theta.
\]
Substituting these expressions into definitions of $\Lambda_1$ and $W_i$ from Lemma \ref{3Snotation}, we obtain (\ref{Lmbda12}), (\ref{WW1}) and (\ref{WW2}).
Substituting (\ref{Lambda1short}) into definition of $P_1$ we prove (\ref{PP1}).\hspace*{\fill}$\Box$


Next lemma deals with auxiliary measures from Theorem \ref{T1}.

 \begin{lemma}\label{3Mlema} Let condition (\ref{ab}) hold. Then, for all $n=1,2,\dots,$
 \begin{eqnarray}
 H-I\!&=&\!(1\!-\!2\A)U+(1\!-\!2\A)0.154U*\Theta,\label{3a0}\\
 H-I\!&=&\!\frac{1}{1-2\A}U+\frac{2\A}{(1-2\A)^2}U^{*2}+ C\A^2U^{*3}*\Theta,\label{3a00}\\
  K\!&=&\!I+C\A(I_1\!-\!I)^{*2}*\Theta.\label{3a4}
 \end{eqnarray}
 \end{lemma}

\textbf{Proof.} For the proof of (\ref{3a0}) observe that
\[H\!-\!I=(1\!-\!2\A)\sum_{j=0}^\infty(2\A)^j(L^{*(j+1)}\!-\!I)=(1\!-\!2\A)(L\!-\!I)+(1\!-\!2\A)(L\!-\!I)*\sum_{j=1}^\infty(2\A)^j\sum_{m=0}^j L^{*m}.\]
We have $L^{*m}\in\F$. Therefore, $\tvnorm{L^{*m}}=1$ and
\[\bigg\|\sum_{j=1}^\infty(2\A)^j\sum_{m=0}^jL^{*m}\bigg\|_{_{TV}}\leq \sum_{j=1}^\infty(2\A)^j\sum_{m=0}^j\tvnorm{L^{*m}}=\sum_{j=1}^\infty(2\A)^j(j\!+\!1)=\frac{4\A-4\A^2}{(1-2\A)^2}\leq 0.154,\]
since $2\A\leq 1/15$ by condition (\ref{ab}).

For the proof of (\ref{3a00}) we note that Fourier transform of $\w H(t)$ allows expression
\[\w H(t)=\frac{(1\!-\!2\A)\w L(t)}{1\!-\!2\A\w L(t)}=1+\frac{\w U(t)}{1\!-\!2\A}\sum_{j=0}^\infty\bigg(\frac{2\A}{1\!-\!2\A}\bigg)^j\w U^j(t).\]
Therefore, due to relations between Fourier transforms and measures,
\begin{eqnarray*}
H&=&I+\frac{U}{1\!-\!2\A}\sum_{j=0}^\infty\bigg(\frac{2\A}{1\!-\!2\A}\bigg)^j U^{*j}\\
&=&I+\frac{1}{1-2\A}U+\frac{2\A}{(1-2\A)^2}U^{*2}+\frac{(2\A)^2}{(1\!-\!2\A)^3}U^{*3}*\sum_{j=0}^\infty\bigg(\frac{2\A}{1\!-\!2\A}\bigg)^j U^{*j}.
\end{eqnarray*}
The proof of (\ref{3a00}) now follows from assumption (\ref{ab}) and simple estimate
\[\frac{4}{(1\!-\!2\A)^3}\bigg\|\sum_{j=0}^\infty\bigg(\frac{2\A}{1\!-\!2\A}\bigg)^j U^{*j}\bigg\|_{_{TV}}\leq\frac{4\cdot 15^3}{14^3}\sum_{j=0}^\infty\frac{1}{7^j}\leq 6.\]

The proof of (\ref{3a4}) is very similar to the proof of (\ref{3a0}). We resp. use Fourier transform and  estimate 
\[\w K(t)=\sum_{j=0}^\infty\bigg(\frac{2\A}{1\!-\!2\A\!+\!2\B}\bigg)^j\w U^j(t),\quad
\sum_{j=1}^\infty\bigg(\frac{2\A}{1\!-\!2\A\!+\!2\B}\bigg)^j\tvnorm{U}^{j-1}\leq C.\]
 \hspace*{\fill}$\Box$


In \cite{SC16}, Lemma 9 the following estimates for closeness of $\Lambda_1$ and $G$ were obtained.
Let
\begin{equation}\label{A0}
A_0=\frac{-2\B^2(1-2\A)}{(1-2\A+2\B)^2}\bigg( (1+2\A)\dirac+\frac{2(1-2\A)K}{1-2\A+2\B} \bigg)*(H-\dirac)^{*2}.
\end{equation}

\begin{lemma} \label{L5} Let condition (\ref{ab}) be satisfied. Then, for all $n=1,2,\dots,$
\begin{eqnarray*}
\lefteqn{\tvnorm{\Lambda_1^{*n}-G^{*n}}\leq C \min\bigg\{\frac{1}{n},n\B^2\bigg\},\qquad
\tvnorm{\Lambda_1^{*n}-G^{*n}*(\dirac\!+\!nA_0)}\leq C \min\bigg\{\frac{1}{n^2},n\B^3\bigg\},}\hskip 15cm\label{lg1}
\\
\lefteqn{\lnorm{\Lambda_1^{*n}-G^{*n}}\leq C \min\bigg\{\frac{1}{n\sqrt{n\B}},n\B^2\bigg\},\quad\lnorm{\Lambda_1^{*n}-G^{*n}*(\dirac\!+\!nA_0)}\leq
 C \min\bigg\{\frac{1}{n^2\sqrt{n\B}},n\B^3\bigg\},
}\hskip 15cm\label{lg2}
\\
\lefteqn{
\wnorm{\Lambda_1^{*n}-G^{*n}}\leq C \min\bigg\{\sqrt{\frac{\B}{n}},n\B^2\bigg\},
\quad\wnorm{\Lambda_1^{*n}-G^{*n}(\dirac\!+\!nA_0)}\leq C \min\bigg\{\frac{\sqrt{\B}}{n^{3/2}},n\B^3\bigg\}.
}\hskip 15cm\label{lg3}
\end{eqnarray*}
\end{lemma}


Next we explore closeness of $D$ and $G$.

\begin{lemma}\label{DG} Let condition (\ref{ab}) hold. Then, for all $n=1,2,\dots,$
\begin{eqnarray}
G^{*n}\!&=&\!C\exp\{0.5n\B(I_1\!-\!I)\}*\Theta,\label{3a1}\\
D^{*n}\!&=&\!\exp\{0.5n\B(I_1\!-\!I)\}*\Theta,\label{3a2}\\
G^{*n}\!&=&\!D^{*n}+Cn\B\A(I_1\!-\!I)^{*4}*\exp\{0.5n\B(I_1\!-\!I)\}*\Theta,\label{3a3}\\
D\!&=&\!I+2\B U/(1\!-\!2\A\!+\!2\B)+C\B^2U^{*2}*\Theta,\label{D1m}\\
G\!&=&\!D*(I+4\A\B U^{*2}/(1\!-\!2\A\!+\!2\B)(1\!-\!2\A)+C\A^2\B U^{*3}*\Theta),\label{3a5}\\
G^{*n}\!&=&\!D^{*n}*(I+4n\A\B U^{*2}/(1\!-\!2\A\!+\!2\B)(1\!-\!2\A))\nonumber\\
&+&\!C\A^2(I_1\!-\!I)^{*4}\exp\{0.2n\B(I_1-I)\}*\Theta.\label{3a6}
\end{eqnarray}
\end{lemma}

\textbf{Proof.} For the proof of (\ref{3a1}) observe that
\[U=\frac{1}{2}(I_1\!-\!I+I_{-1}\!-\!I)=\frac{1}{4}(I_1\!-\!I)^{*2}*I_{-1}+\frac{1}{4}(I_{-1}\!-\!I)^{*2}*I_1.\]
Substituting these expressions into (\ref{3a0}) we get
\begin{eqnarray*}
G^{*n}\!&=&\!\exp\{\omega(U+0.154U*\Theta)\}\\
\!&=&\!\exp\Big\{\frac{\omega}{2}(I_1\!-\!I)+\frac{0.154\omega}{4}(I_1\!-\!I)^{*2}*\Theta\Big\}*\exp\Big\{\frac{\omega}{2}(I_{-1}\!-\!I)+\frac{0.154\omega}{4}(I_{-1}\!-\!I)^{*4}*\Theta\Big\}.
\end{eqnarray*}
Here
\[\omega=\frac{2n\B(1-2\A)^2}{1-2\A+2\B}.\]
Then by Lemma \ref{L2.5}
\[G^{*n}=C\exp\{0.7305\,\omega(I_1-I)\}*\exp\{0.7305\,\omega(I_{-1}-I)\}*\Theta=C\exp\{0.7305\,\omega(I_1-I)\}*\Theta .\]
Observing that, due to condition (\ref{ab}),
$\omega\geqslant 0.5n\B$
 we complete the proof of (\ref{3a1}). Similarly, (\ref{3a2}) follows from
$n\B/(1-2\A+2\B)>0.5n\B$.

For the proof of (\ref{3a3}) we apply Lemma \ref{bazine} with
\[\bar{M}=\frac{2n\B(1-2\A)}{1-2\A+2\B}(H-I),\qquad \bar{V}=\frac{2n\B}{1-2\A+2\B}U\]
arriving at
\begin{equation}\label{GtDn}
G^{*n}- D^{*n}=\eec^{\bar{M}}-\eec^{\bar{V}}=(\bar{M}-\bar{V})*\int_0^1\eec^{\tau \bar{M}+(1-\tau)\bar{V}}\dd\tau.
\end{equation}
From (\ref{3a1}) and (\ref{3a2}) it follows that
\[\int_0^1\eec^{\tau \bar{M}+(1-\tau)\bar{V}}\dd\tau=C\eec^{0.5n\B(I_1-I)}*\Theta.\]
From (\ref{3a00}) and (\ref{X2}) we get
\begin{eqnarray*}
\bar{V}-\bar{M}\!&=&\!\frac{2n\B}{1-2\A+2\B}((1-2\A)(H-I)-U)=
\frac{4n\A\B}{(1-2\A+2\B)(1-2\A)}U^{*2}*\Theta\\
\!&=&\!Cn\A\B(I_1-I)^{*4}*\Theta.
\end{eqnarray*}
Substituting the last two equalities into (\ref{GtDn}) we complete the proof of (\ref{3a3}). Estimate (\ref{D1m}) follows from Lemma \ref{bazine}.
Similarly, from (\ref{3a00}) and Lemma \ref{bazine} it follows that
\begin{eqnarray*}
G\!&=&\!\exp\Big\{\frac{2\B(1\!-\!2\A)}{1\!-\!2\A\!+\!2\B}\Big(\frac{1}{1-2\A}U+\frac{2\A}{(1-2\A)^2}U^{*2}+ C\A^2U^{*3}*\Theta\Big)\Big\}\\
\!&=&\!D*\exp\Big\{\frac{4\A\B U^{*2}}{(1\!-\!2\A\!+\!2\B)(1\!-\!2\A)}+C\A^2\B U^{*3}*\Theta\Big\}\\
\!&=&\!D*(I+4\A\B u^{*2}/(1\!-\!2\A\!+\!2\B)(1\!-\!2\A)+C\A^2\B U^{*3}*\Theta),
\end{eqnarray*}
which completes the proof of (\ref{3a5}).

For the proof of (\ref{3a6}) observe that from (\ref{3a5}) it follows that
\begin{equation}\label{GD}
G-D=C\A\B U^{*2}*\Theta= C\A\B(I_1-I)^{*4}*\Theta.
\end{equation}
Therefore, from Lemma \ref{bergidentL}, (\ref{3a1}), (\ref{3a2}) and (\ref{C.3}) it follows that
\begin{eqnarray*}
G^{*n}&=&D^{*n}+nD^{*(n-1)}*(G-D)+Cn^2\A^2\B^2(I_1-I)^{*8}\exp\{0.5n\B(I_1-I)\}*\Theta\\
&=&D^{*n}+nD^{*(n-1)}*(G-D)+C\A^2(I_1-I)^{*4}*\exp\{0.2n\B(I_1-I)\}*\Theta.
\end{eqnarray*}
From (\ref{3a5})  we obtain
\[nD^{*(n-1)}*(G-D)=D^{*n}4n\A\B U^{*2}/(1\!-\!2\A\!-\!2\B)(1\!-\!2\A)+Cn\A^2\B(I_1-I)^{*6}*\exp\{0.5n\B(I_1-I)\}*\Theta.\]
By (\ref{C.2})
\[Cn\A^2\B(I_1-I)^{*6}*\exp\{0.5n\B(I_1-I)\}*\Theta=C\A^2(I_1-I)^{*4}*\exp\{0.2n\B(I_1-I)\}*\Theta.\]
Collecting the last three expressions we complete the proof of (\ref{3a6}). \hspace*{\fill}$\Box$ \\


Next lemma is needed for asymptotic expansions. Let
\begin{equation}\label{delta}
\delta:=-\frac{2\B^2}{(1\!-\!2\A\!+\!2\B)^2}\bigg(\frac{1\!+\!2\A}{1\!-\!2\A}+\frac{2}{1\!-\!2\A\!+\!2\B}\bigg).
\end{equation}

\begin{lemma}\label{GDA} Let condition (\ref{ab}) hold. Then, for all $n=1,2,\dots,$
\begin{eqnarray*}
nG^{*n}*A_0&=&n\delta G^{*n}*U^{*2}+C(\A\B\!+\!\B^2)(I_1-I)^{*4}*\exp\{0.2n\B(I_1-I)\}*\Theta\\
&=&n\delta D^{*n}*U^{*2}+C(\A\B\!+\!\B^2)(I_1-I)^{*4}*\exp\{0.2n\B(I_1-I)\}*\Theta,\\
G^{*n}*(I\!+\!nA_0)&=&D^{*n}*(I+nA_2)+C(\A\B\!+\!\B^2)(I_1-I)^{*4}*\exp\{0.2n\B(I_1-I)\}*\Theta.
\end{eqnarray*}
\end{lemma}

\textbf{Proof.} From (\ref{3a0}), (\ref{3a4}), (\ref{3a1}) and (\ref{C.2}) we have
\begin{eqnarray*}
nG^{*n}*A_0-n\delta G^{*n}*U^{*2}&=&Cn\A\B^2(I_1-I)^{*6}*\exp\{0.5n\B(I_1-I)\}*\Theta\\
&=&C\A\B(I_1-I)^{*4}*\exp\{0.2n\B(I_1-I)\}*\Theta.
\end{eqnarray*}
Remaining expressions are proved by using (\ref{3a5}) and (C.3).\hspace*{\fill}$\Box$

\section{Proofs}

\textbf{Proof of Theorem \ref{T1GC}}. We have
\begin{eqnarray}\label{prT1a}
F_n-D^{*n}&=&P_1*(\Lambda_1^{*n}-G^{*n})*W_1+P_1*(G^{*n}-D^{*n})*W_1+(P_1-I)*D^{*n}*W_1\nonumber\\
&+&D^{*n}*(W_1-I)+P_1*\Lambda_2^{*n}*W_2.
\end{eqnarray}

From Lemma \ref{L5} and (\ref{nrmrel}) it follows that
\[\wnorm{P_1*(\Lambda_1^{*n}-G^{*n})*W_1}\leq C\wnorm{\Lambda_1^{*n}-G^{*n}}\leq C\B^{1/2}n^{-1/2}.\]
Estimates (\ref{nrmrel}), (\ref{locWToTV}), (\ref{C.3}) and Lemma (\ref{3a3}) allow to obtain
\begin{eqnarray*}
\wnorm{P_1*(G^{*n}\!-\!D^{*n})*W_1}\!&\leq&\!
  C\wnorm{G^{*n}\!-\!D^{*n}}\leq Cn\A\B\wnorm{(I_1\!-\!I)^{*4}*\exp\{0.5n\B(I_1\!-\!I)\}*\Theta}\\
  &\leq& Cn\A\B\tvnorm{(I_1\!-\!I)^{*3}*\exp\{0.5n\B(I_1\!-\!I)\}}\\
  &\leq& C\A n^{-1/2}\B^{-1/2}.
\end{eqnarray*}
Similarly, from (\ref{PP1}) and (\ref{C.2}) we get
\begin{eqnarray*}
\wnorm{(P_1-I)*D^{*n}*W_1}\!&\leq&\! C(\A+\B)\wnorm{(I_1\!-\!I)^{*2}*D^{*n}}=C(\A+\B)\tvnorm{(I_1\!-\!I)*D^{*n}}\\
\!&\leq&\!C(\A+\B)\tvnorm{(I_1\!-\!I)**\exp\{0.5n\B(I_1\!-\!I)\}}\leq C(\A+\B)n^{-1/2}\B^{-1/2} \end{eqnarray*}
and from (\ref{WW1}) and (\ref{C.2})
\[\wnorm{D^{*n}*(W_1-I)}\leq C(\A+\B)n^{-1/2}\B^{-1/2}.\]
Next, from Lemma \ref{3Snotation}
\begin{eqnarray*} 
\wnorm{P_2*\Lambda_2^{*n}*W_2}\!&\leq&\! C(\A+\B)\wnorm{(I_1\!-\!I)^{*2}*\Lambda_2^{*n}}=C(\A+\B)\tvnorm{(I_1\!-\!I)*\Lambda_2^{*n}}\\
\!&\leq&\!C(\A+\B)\tvnorm{\Lambda_2}^n\leq C(\A+\B)0.2^n\leq C(\A+\B)n^{-1/2}.\end{eqnarray*}
Collecting all above estimates we get
\[\wnorm{F_n-D^{*n}}\leq C\sqrt{\frac{\B}{n}}\bigg(1+\frac{\A}{\B}\bigg).\]
To complete the proof for Wasserstein metric observe that
\[1+\frac{\A}{\B}=2+\frac{\A-\B}{\B}\leq 2\bigg(1+\frac{\ab{\A-\B}}{\B}\bigg).\]
Proofs for total variation and point metrics are very similar and, therefore, omitted. \hspace*{\fill}$\Box$\\


\textbf{Proof of Theorem \ref{T2}.} We can write identity
\[ F_n-D^{*n}*(I\!+\!A_1\!+\!nA_2)=J_1+J_2+J_3+J_4+J_5+J_6,\]
\begin{eqnarray*}
J_1\!&=&\!P_1*(\Lambda_1^{*n}-G^{*n}*(I\!+\!nA_0))*W_1,\\
J_2\!&=&\!P_1*(G^{*n}*(I\!+\!nA_0)-D^{*n}*(I\!+\!nA_2))*W_1,\\
J_3\!&=&\!(P_1-I-2p_2(\B-\A)U/(1\!-\!2\A\!+\!2\B))*D^{*n}*(I\!+\!nA_2)*W_1,\\
J_4\!&=&\!\bigg(I+\frac{2p_2(\B\!-\!\A)U}{1\!-\!2\A\!+\!2\B}\bigg)*D^{*n}*(I\!+\!nA_2)*\bigg(W_1-I-\frac{2(\A\!-\!\B)(1-2\A)U}{(1\!-\!2\A\!+\!2\B)^2}\bigg),\\
J_5\!&=&\!D^{*n}*\bigg(\bigg(I+\frac{2p_2(\B\!-\!\A)U}{1\!-\!2\A\!+\!2\B}\bigg)*\bigg(I+\frac{2(\A\!-\!\B)(1-2\A)U}{(1\!-\!2\A\!+\!2\B)^2}\bigg)
*(I\!+\!nA_2)-(I\!+\!A_1\!+\!nA_2)\bigg),\\
J_6\!&=&\!P_2*\Lambda_2^{*n}*W_2.
\end{eqnarray*}
The proof is now very similar to the proof of Theorem \ref{T1GC}. From Lemma \ref{L5} and (\ref{baznrm}) we get
$\lnorm{J_1}\leq Cn^{-5/2}\B^{-1/2}$. Applying Lemma \ref{GDA}, (\ref{locWToTV}), (\ref{3a2})  and (\ref{C.3}) we obtain
\begin{eqnarray*}\lnorm{J_2}\!&\leq&\! C(\A\B+\B^2)\lnorm{(I_1\!-\!I)^{*4}*\exp\{0.2n\B(I_1\!-\!I)\}}\\
\!&\leq&\! C(\A\B+\B^2)\tvnorm{(I_1\!-\!I)^{*5}*\exp\{0.2n\B(I_1\!-\!I)\}}\\
\!&\leq&\!
Cn^{-2}(n\B)^{-1/2}(1+\A/\B)\leq Cn^{-2}(n\B)^{-1/2}(2+2\A^2/\B^2).
\end{eqnarray*}
Observe that from (\ref{X2}), (\ref{3a2}) and (\ref{C.2}) it follows that
\begin{eqnarray}\label{nDA2} 
D^{*n}*nA_2&=& Cn(\B^2+\A\B)(I_1-I)^{*4}*\exp\{0.5n\B(I_1\!-\!I)\}*\Theta\nonumber\\
&=&C(\A+\B)(I_1\!-\!I)^{*2}*\exp\{0.2n\B(I_1\!-\!I)\}*\Theta\end{eqnarray}
and
\[D^{*n}*(I\!+\!nA_2)=C\exp\{0.2n\B(I_1\!-\!I)\}*\Theta.\]
Therefore applying (\ref{PP1}), (\ref{WW1}), (\ref{3a2}) and (\ref{C.3}) we  prove that
\[
\lnorm{J_3}\leq  C(\A\B+\B^2)\lnorm{(I_1\!-\!I)^{*4}*\exp\{0.2n\B(I_1\!-\!I)\}}
\leq Cn^{-2}(n\B)^{-1/2}(2+2\A^2/\B^2)
\]
and
\[\lnorm{J_4}\leq C(\A^2+\B^2)\lnorm{(I_1\!-\!I)^{*4}*\exp\{0.2n\B(I_1\!-\!I)\}}
\leq Cn^{-2}(n\B)^{-1/2}(1+\A^2/\B^2).\]
Expression (\ref{nDA2}) allows to write
\[J_5=C(\A+\B)^2(I_1-I)^{*4}*\exp\{0.2n\B(I_1\!-\!I)\}*\Theta.\]
Therefore, by (\ref{locWToTV})  and (\ref{C.3})
\[\lnorm{J_5}\leq C(\A+\B)^2\tvnorm{(I_1-I)^{*5}*\exp\{0.2n\B(I_1\!-\!I)\}}\leq Cn^{-2}(n\B)^{-1/2}(1+\A^2/\B^2).\]
Finally, from  Lemma \ref{3Snotation}
\begin{eqnarray*} 
\lefteqn{\lnorm{P_2*\Lambda_2^{*n}*W_2}\leq C(\A+\B)\tvnorm{(I_1\!-\!I)*\Lambda_2^{*n}}
\leq C(\A+\B)0.2^n}\hskip 1cm\\
\!&\leq&\! C(\A+\B)n^{-5/2}\leq Cn^{-5/2}\B^{-1/2}(1+\A/\B)\leq Cn^{-5/2}\B^{-1/2}(1+\A^2/\B^2).\end{eqnarray*}
Collecting all above estimates we get
\[\lnorm{F_n-D^{*n}*(I\!+\!A_1\!+\!nA_2)}\leq Cn^{-2}(n\B)^{-1/2}(1+\A^2/\B^2),\]
which is equivalent to Theorem's statement for point metric, since $2(\A-\B)\B\leq (\A-\B)^2+\B^2$ and, therefore
\[1+\frac{\A^2}{\B^2}=2+\frac{(\A\!-\!\B)^2+2(\A\!-\!\B)\B}{\B^2}\leq 3+\frac{2(\A\!-\!\B)^2}{\B^2}\leq 3\bigg(1+\frac{(\A\!-\!\B)^2}{\B^2}\bigg).\]
Proof for total variation and Wasserstein norms is very similar and, therefore, omitted. \hspace*{\fill}$\Box$\\


\textbf{Proof of Theorem \ref{T1star}}. We have
\[\tvnorm{F_n-D^{*n}}\leq\tvnorm{F_n-D^{*n}*(I\!+\!A_1\!+\!nA_2)}+\tvnorm{A_1*D^{*n}}+n\tvnorm{A_2*D^{*n}}.\]
We  can apply Theorem \ref{T2} to estimate the first summand. 
For the other estimates we use Lemma \ref{Roslema}. Recall that $D=\exp\{\l(I_1\!-\!I)+\l(I_{-1}\!-\!I)\}$, $\l=\B/(1-2\A+2\B)$. Observe that $U=-(I_1\!-\!I)*(I_{-1}\!-\!I)/2$. Therefore by (C.2)
\begin{eqnarray}\label{k1}
\tvnorm{A_1*D^{*n}}\!&\leq&\!\frac{\ab{\A\!-\!\B}}{1\!-\!2\A\!+\!2\B}\max\bigg(\frac{1\!-\!2\A}{1\!-\!2\A\!+\!2\B},p_2\bigg)\tvnorm{(I_1\!-\!I)*\exp\{n\l(I_1\!-\!I)\}}\nonumber\\
&\times&\!\tvnorm{(I_{-1}\!-\!I)*\exp\{n\l(I_1\!-\!I)\}}\leq 2\ab{\A\!-\!\B}(\eec n\B)^{-1}.
\end{eqnarray}
Similarly
\begin{eqnarray}\label{k2}
n\tvnorm{A_2*D^{*n}}\!&\leq&\!\frac{2\B n}{(1\!-\!2\A\!+\!2\B)^2}(2\ab{\A\!-\!\B}+\B)\cdot\frac{1}{4}\tvnorm{(I_1\!-\!I)^{*2}*\exp\{n\l(I_1\!-\!I)\}}\nonumber\\
\!&\times&\!\tvnorm{(I_{-1}\!-\!I)^{*2}*\exp\{n\l(I_{-1}\!-\!I)\}}\nonumber\\
\!&\leq&\!\frac{9}{\eec^2 n}\cdot\frac{\ab{\A\!-\!\B}}{\B}+\frac{9}{2\eec^2 n}.
\end{eqnarray}
From (\ref{k1}) and (\ref{k2}) theorem's statement for total variation norm follows. 

For local norm observe that by (\ref{ab}) $1\!-\!2\A\!+\!2\B\leq 16/15$ and by Lemma \ref{zloc}
\begin{eqnarray*}\label{k3}
\lnorm{A_1*D^{*n}}\!&\leq&\!\frac{2\ab{\A\!-\!\B}}{1\!-\!2\A\!+\!2\B}\lnorm{U*\exp\{2n\l U\}}\leq \frac{4\ab{\A\!-\!\B}}{n\B\sqrt{n\B}}\bigg(\frac{3/2}{2\eec}\bigg)^{3/2}\sqrt{\frac{16}{15}}\leq
\frac{0.6\ab{\A\!-\!\B}}{n\B\sqrt{n\B}}
\end{eqnarray*}
and 
\begin{eqnarray*}\label{k4}
n\lnorm{A_2*D^{*n}}\!&\leq&\!\frac{2\B n}{(1\!-\!2\A\!+\!2\B)^2}(2\ab{\A\!-\!\B}+\B)\lnorm{U^{*2}*\exp\{2n\l U)\}}\\
\!&\leq&\!\frac{4 (2\ab{\A\!-\!\B}+\B)}{n\B\sqrt{n\B}} \bigg(\frac{5}{4\eec}\bigg)^{5/2}\sqrt{\frac{16}{15}}\leq
\frac{0.6 (2\ab{\A\!-\!\B}+\B)}{n\B\sqrt{n\B}}.
\end{eqnarray*}

For Wasserstein norm we note that $1\!-\!2\A\!+\!2\B\geq 14/15$ and use (\ref{locWToTV}), (\ref{C.2}) and (\ref{C.3})
\begin{eqnarray*}
\wnorm{A_1*D^{*n}}\!&\leq&\!\frac{\ab{\A\!-\!\B}}{1\!-\!2\A\!+\!2\B}\wnorm{(I_1\!-\!I)*\exp\{n\l (I_1\!-\!I)\}*(I_1\!-\!I)*\exp\{n\l(I_1\!-\!I)\}}\\
\!&\leq&\!\frac{\ab{\A\!-\!\B}}{1\!-\!2\A\!+\!2\B}\tvnorm{(I_1\!-\!I)*\exp\{n\l(I_1\!-\!I)\}}\\
\!&\leq&\!\frac{\ab{\A\!-\!\B}}{\sqrt{n\B}}\sqrt{\frac{15}{14}\cdot\frac{2}{\eec}}\leq \frac{0.888\ab{\A\!-\!\B}}{\sqrt{n\B}}
\end{eqnarray*}
and
\begin{eqnarray*}
n\wnorm{A_2*D^{*n}}\!&\leq&\!\frac{\B n}{2(1\!-\!2\A\!+\!2\B)^2}(2\ab{\A\!-\!\B}+\B)\tvnorm{(I_1\!-\!I)*\exp\{n\l(I_1\!-\!I)\}}\nonumber\\
\!&\times&\!\tvnorm{(I_{-1}\!-\!I)^{*2}*\exp\{n\l(I_{-1}\!-\!I)\}}\nonumber\\
\!&\leq&\!\frac{3}{\eec\sqrt{2\eec}}\sqrt{\frac{15}{14}}\,\frac{(2\ab{\A\!-\!\B}+\B)}{\sqrt{n\B}}\leq \frac{0.49(2\ab{\A\!-\!\B}+\B)}{\sqrt{n\B}}.
\end{eqnarray*}
 Combining the above estimates with Theorem \ref{T2} we complete the proof of Theorem \ref{T1star}.
\hspace*{\fill}$\Box$\\


\textbf{Proof of Theorem \ref{lrLr}}.  Obviously, $\tvnorm{M}=\norm{M}_1$, $\wnorm{M}=\ab{M}_1$.
For $r>1$ we have
\[\dlrnorm{F_n- D^{*n}}\leq (\lnorm{F_n- D^{*n}})^{(r-1)/r}(\tvnorm{F_n- D^{*n}})^{1/r}.\]
Observe that
\[\sup_k\ab{F_n\{(-\infty,k]\}-D^{*n}\{(-\infty,k]\}}\leq\sup_A\ab{F_n\{A\}-D^{*n}\{A\}}\leq\tvnorm{F_n-D^{*n}},\]
where supremum is taken by all Borel sets.
Therefore,
 
\[\Lnorm{F_n- D^{*n}}\leqslant(\tvnorm{F_n- D^{*n}})^{(r-1)/r}(\wnorm{F_n- D^{*n}})^{1/r}.\]

It remains to use estimates from Theorem \ref{T1GC}. \hspace*{\fill}$\Box$

     {\small

\end{document}